%% file: Gehonico.tex
\input amstex
\documentstyle{amsppt}
\TagsOnLeft
\topmatter
\title Geometric and homotopy theoretic methods in Nielsen coincidence theory \endtitle
\rightheadtext{Methods in Nielsen coincidence theory}
\author Ulrich Koschorke\endauthor
\leftheadtext{Ulrich Koschorke}
\address Universit\"at Siegen,
Emmy Noether Campus, Walter-Flex-Str. 3,
D-57068 Siegen, Germany
\endaddress
\email koschorke\@mathematik.uni-siegen.de \endemail

\vskip2in
\abstract \ In classical fixed point and coincidence theory the notion of Nielsen numbers has proved to be extremely fruitful. Here we extend it to pairs \ $(f_1, f_2)$ \ of maps between manifolds of arbitrary dimensions. This leads to estimates of the minimum numbers \ $MCC (f_1, f_2)$ \ (and \ $MC (f_1, f_2)$, resp.) of pathcomponents (and of points, resp.) in the coincidence sets of those pairs of maps which are homotopic to \ $(f_1, f_2)$.  Furthermore we deduce finiteness conditions for \ $MC (f_1, f_2)$. As an application we compute both minimum numbers explicitly in four concrete geometric sample situations.

The Nielsen decomposition of a coincidence set is induced by the decomposition of a certain path space \ $E (f_1, f_2)$ \ into pathcomponents. Its higher dimensional topology captures further crucial geometric coincidence data.

An analoguous approach can be used to define also Nielsen numbers of certain link maps.
\endabstract
\endtopmatter


\input xy
\input xymatrix
\input xyarrow
\input xycurve
\input graphicx

\define\coker{\operatorname{coker}}
\define\incl{\operatorname{incl}}

\define\coll{\operatorname{coll}}



\document

\input gehoni01.tex

\input gehoni02.tex
\input gehoni03.tex
\input gehoni04.tex
\input gehoni05.tex

\input gehoni06.tex
\input gehonref.tex

\enddocument

%% file: gehoni01.tex
\specialhead 1.\ \ Introduction and discussion of results
\endspecialhead

Throughout this paper \ $f_1, f_2 : M \to N$ \ denote two (continuous) maps between the smooth connected manifolds \ $M$ \ and \ $N$ \ without boundary, of strictly positive dimensions \ $m$ \ and \ $n$, \ resp., \ $M$ \ being compact.

We would like to measure how small (or simple in some sense) the coincidence locus
$$
C (f_1, f_2) \ := \ \{ x \in M \ | \ f_1 (x) = f_2 (x) \}
\tag 1.1
$$
can be made by deforming \ $f_1$ \ and \ $f_2$ \ via homotopies. Classically one considers the \ {\bf m}inimum number of \ {\bf c}oincidence points
$$
MC (f_1, f_2) \ := \ \min \{ \# C (f'_1, f'_2) | f'_1 \sim f_1, f'_2 \sim f_2\}
\tag 1.2
$$
(cf.\ \cite{BGZ}, 1.1). It coincides with the minimum number \ $\min \{ \# C (f'_1, f_2) | f'_1 \sim f_1\}$ \ where only \ $f_1$ \ is modified by a homotopy (cf.\ \cite{Br}). In particular, in topological fixed point theory (where \ $M = N$ \ and \ $f_2$ \ is the identity map) this minimum number is the principal object of study (cf.  \cite{B}, p.\ 9).

In higher codimensions, however, the coincidence locus is generically a manifold of dimension \ $m - n > 0$, \ and \ $MC (f_1, f_2)$ \ is often infinite (see e.g.\ examples I and III  below). Thus it seems more meaningful to study the \ {\bf m}inimum number of \ {\bf c}oincidence \ {\bf c}omponents
$$
MCC (f_1, f_2) \ := \ \min \{ \# \pi_0 (C (f'_1, f'_2)) | f'_1 \sim f_1, f'_2 \sim f_2\}
\tag 1.3
$$
where \ $\# \pi_0 (C (f'_1, f'_2))$ \ denotes the (generically finite) number of \ {\it pathcomponents} \ of the indicated coincidence subspace of \ $M$.

\subhead Question \endsubhead \ How big are \ $MCC (f_1, f_2)$ \ and \ $MC (f_1, f_2)$? \ In particular, when do these invariants vanish, i.e.\ when can the maps \ $f_1$ \ and \ $f_2$ \ be deformed away from one another?
\vskip2mm

In this paper we discuss lower bounds for \ $MCC (f_1, f_2)$ \ and  geometric obstructions to \ $MC (f_1, f_2)$ \ being trivial or finite.

A careful investigation of the differential topology of generic coincidence submanifolds yields the normal bordism classes (cf.\ 4.6 and 4.7)
$$
\omega (f_1, f_2) \ \in \ \Omega_{m -n} (M; \varphi)
\tag 1.4
$$
and
$$
\widetilde\omega (f_1, f_2) \ \in \ \Omega_{m -n} (E (f_1, f_2); \ \widetilde\varphi)
\tag 1.5
$$
as well as a sharper (\lq\lq nonstabilized\rq\rq) version
$$
\omega^\# (f_1, f_2) \ \in \ \Omega^\# (f_1, f_2)
\tag 1.6
$$
of \ $\widetilde\omega (f_1, f_2)$ \ (cf.\ 4.10). Here the pathspace
$$
E (f_1, f_2) \ := \ \{ (x, \theta) \in M \times N^I \ | \ \theta (0) = f_1 (x), \theta (1) = f_2 (x) \}
\tag 1.7
$$
(compare \S\ 2), also known as (a kind of) \ {\it homotopy equalizer of \ $f_1$ \ and \ $f_2$}, \ plays a crucial role. In general it has a very rich topology involving both \ $M$ \ and the loop space of \ $N$. \ Already the set \ $\pi_0 (E (f_1, f_2))$ \ of pathcomponents can be huge -- it corresponds bijectively to the Reidemeister set
$$
R (f_1, f_2) \ \ = \ \ \pi_1 (N) / \text{Reidemeister equivalence}
\tag 1.8
$$
(compare \cite{BGZ}, 3.1 and our proposition 2.3 below) which is of central importance in classical Nielsen theory. Thus it is only natural to define a \lq\lq Nielsen number\rq\rq\ \ $N (f_1, f_2)$ \ (and a sharper version \ $N^\# (f_1, f_2)$, resp.) to be the number of those (\lq\lq essential\rq\rq ) pathcomponents which contribute nontrivially to the bordism class \ $\widetilde\omega (f_1, f_2)$ \ (and to \ $\omega^\# (f_1, f_2)$, resp.), cf.\ 4.9 and 4.10 below.

\proclaim{Theorem 1.9} \ {\rm(i)} \ The integers \ $N (f_1, f_2)$ \ and \ $N^\# (f_1, f_2)$ \ depend only on the homotopy classes of \ $f_1$ \ and \ $f_2$;

{\rm (ii)} \ \ $N (f_1, f_2) = N (f_2, f_1)$ \ and \ $N^\# (f_1, f_2) = N^\# (f_2, f_1)$;

{\rm (iii)} \ $0 \le N (f_1, f_2) \ \le \ N^\# (f_1, f_2) \ \le \ MCC (f_1, f_2) \ \le MC (f_1, f_2);$ \newline
if \ $n \ne 2$, \ then also \ $MCC (f_1, f_2) \ \le \ \# R (f_1, f_2)$;

{\rm(iv)} \ \ if \ $m = n$ \ then \ $N (f_1, f_2) = N^\# (f_1, f_2)$ \ coincides with the classical Nielsen number {\rm (}cf.\ \cite{BGZ}, definition {\rm3.6)}.
\endproclaim

\remark{Remark 1.10} \ In various situations some of the estimates spelled out in part (iii) of this theorem are known to be sharp (compare also \cite{K 3}). E.g.\ in the selfcoincidence setting (where \ $f_1 = f_2$) we have always \ $MCC (f_1, f_2) \le 1$ \ (since here \ $C (f_1, f_2) = M$). In the \lq\lq root setting\rq\rq\ (where \ $f_2$ \ maps to a constant value \ $* \in N$) \ all Nielsen classes are simultaneously essential or inessential (since our \ $\omega$-invariants are always compatible with homotopies of \ $(f_1, f_2)$ \ and hence, in this particular case, with the action of \ $\pi_1 (N, *)$, \ cf.\ the discussion in \cite{K 3} following 1.11). Therefore in both settings \ $MCC (f_1, f_2)$ \ is equal to the Nielsen number \ $N (f_1, f_2)$ \ provided \ $\widetilde\omega (f_1, f_2) \ne 0$ \ (and \ $n \ne 2$ \ if \ $f_2 \equiv *$).
\endremark

Further geometric and homotopy theoretic considerations allow us to determine the Nielsen and minimum numbers explicitly in several concrete sample situations (for proofs see \S\ 6 below).

\example{Example I} \ Given integers \ $q > 1$ \ and \ $r$, \ let \ $N = \Bbb C P (q)$ \ be $q$-dimensional complex projective space, let \ $M = S (\otimes^r_{\Bbb C} \lambda_{\Bbb C})$ \ be the total space of the unit circle bundle of the \ $r^{th}$ tensor power of the canonical complex line bundle, and let \ $f : M \to N$ \ denote the fiber projection. Then
$$
N (f, f) = N^\# (f, f) = MCC (f, f) = \left\{ \aligned
0 \ \ \ &\text{if} \ \ q \equiv - 1 (r) \ \text{and} \ q \equiv 1 (2), \\
1 \ \ \ &\text{else} \ \ \ \ \ ;
\endaligned \right.
$$
and
$$
MC (f, f) \ \ \ = \ \ \ \left\{ \aligned
0 \ \ \ &\text{if} \ \ q \equiv - 1 (r) \ \text{and} \ q \equiv  1 (2) , \\
1 \ \ \ &\text{if} \ \ q \equiv - 1 (r) \ \text{and} \ q \equiv  0 (2) , \\
\infty \ \ &\text{if} \ \ q \not\equiv - 1 (r) \ \ .
\endaligned \right.
$$  \hfill $\square$

\endexample

As was shown above (cf.\ 1.10), in any selfcoincidence situation (where \ $f_1 = f_2)$ \ $MCC (f_1, f_2)$ \ must be \ $0$ \ or \ $1$ \ and it remains only to decide which value occurs. In the previous example this can be settled by the normal bordism class \ $\omega (f, f) \in \Omega_1 (M; \varphi)$, a weak form of \ $\widetilde\omega (f, f)$ \ which, however, captures a delicate (\lq\lq second order\rq\rq) \ $\Bbb Z_2$--aspect as well as the dual of the classical first order obstruction. Already in this simple case standard methods of singular (co)homology theory yield only a necessary condition for \ $MCC (f_1, f_2)$ \ to vanish (cf.\ \cite{DG}, 2.2). In higher codimensions \ $m - n$ \ the advantage of the normal bordism approach can be truely dramatic.

\example{Example II} \ Given natural numbers \ $k < r$, \ let \ $M = V_{r, k}$ \ (and \ $N = G_{r, k,},$ ß\ resp.) be the Stiefel manifold of orthonormal \ $k$--frames (and the Grassmannian of \ $k$--planes, resp.) in \ $\Bbb R^r$. Let \ $f : M \to N$ \ map a frame to the plane it spans.

Assume \ $r \ge 2k \ge 2$. \ Then
$$
N (f, f) = N^\# (f, f) = MCC (f, f) = MC (f, f) = \left\{
\aligned
0 \ \ \ &\text{if} \ \ \omega (f, f) = 0 \ , \\
1 \ \ \ &\text{else} \ \ \ \ \ .
\endaligned  \right.
$$
Here  the normal bordism obstruction \ $\omega (f, f) \in \Omega_{m -n} (M; \varphi)$ \ (cf.\ 4.7) contains precisely as much information as its \lq\lq highest order component\rq\rq
$$
2 \chi (G_{r, k}) \cdot [SO (k)] \ \ \in \ \ \Omega^{fr}_{m -n} \ \cong \ \pi^S_{m -n}
\tag 1.11
$$
where \ $[SO (k)]$ \ denotes the framed bordism class of the Lie group \ $SO (k)$, \ equipped with a left invariant parallelization; the Euler number \ $\chi (G_{r, k})$ \ is easily calculated: it vanishes if \ $k \not\equiv r \equiv 0 (2)$ \ and equals $ \left( \matrix [r/2] \\[k/2] \endmatrix \right)$ \ otherwise. Without loosing its geometric flavor, our original question translates here -- via the Pontryagin-Thom isomorphism -- into deep problems of homotopy theory (compare the discussion in the introduction of \cite{K 2}). Fortunately powerful methods are available in homotopy theory which imply e.g.\ that \ $MCC (f, f) = MC (f, f) = 0$ \ if \ $k$ \ is even or \ $k = 7$ \ or \ $9$ \ or \ $\chi (G_{r, k}) \equiv 0 (12)$; \ however, if \ $k = 1$ \ and \ $r \equiv 1 (2)$, or  if \ $k = 3$ \ and \ $r \not\equiv 1 (12)$ \ is odd, or if \ $k = 5\ $ and \ $r \equiv 5 (6)$, \ then \ $MCC (f,f) = MC (f, f) = 1$.

These results seem to be entirely out of the reach of the methods of singular (co)homology theory since we would have to deal here with obstructions of order \ $m - n + 1 = \frac{k (k -1)}{2}  + 1$. \hfill $\square$
\endexample

\example{Example III} \ Let \ $N$ be the torus \ $(S^1)^n$ \ and let \ $\iota_1, \dots, \iota_n$ \ denote the canonical generators of \ $H^1 ((S^1)^n; \Bbb Z)$. If the homomorphism
$$
f_{1*} \ - \ f_{2*} \ : \ H_1 (M; \Bbb Z) \ \longrightarrow \ H_1 ((S^1)^n; \ \Bbb Z)
\tag 1.12
$$
has an infinite cokernel (or, equivalently, the rank of its image is strictly smaller than \ $n$), then
$$
N (f_1, f_2) = N^\# (f_1, f_2) = MCC (f_1, f_2) = MC (f_1, f_2) = 0 .
$$
On the other hand, if the cup product
$$
\prod^n_{j = 1} (f^*_1 - f^*_2) (\iota_j) \ \in \ H^n (M; \Bbb Z)
$$
is nontrivial then \ $MC (f_1, f_2) = \infty$ \ whenever \ $m > n$; \ if in addition \ $n \ne 2$, then \ $MCC (f_1, f_2)$ \ equals the (finite) cardinality of the cokernel of \ $f_{1 *} - f_{2 *}$ \ (cf. 1.12).

In the special case when \ $N$ is the unit circle \ $S^1$ \ we have:

\flushpar
$MCC (f_1, f_2) = MC (f_1, f_2) = 0$ \ if \ $f_1$ \ is homotopic to \ $f_2$; otherwise \ $MCC (f_1, f_2)$ \ $ = \# \coker (f_{1 *} - f_{2 *})$, \ but (if $m > 1$) \ \ $MC (f_1, f_2) = \infty$. \hfill $\square$
\endexample

An important special case of our invariants are the degrees
$$
\deg^\# (f) := \omega^\# (f, *), \ \ \ \widetilde\deg (f) := \widetilde\omega (f, *) \ \ \text{and} \ \ \deg (f) := \omega (f, *)
\tag 1.13
$$
of a given map \ $f : M \to N$ \ (here \ $*$ \ denotes a constant map).

\example{Example IV (homotopy groups)} \ Let \ $M$ \ be the sphere \ $S^m$; \ in view of the previous example we may also assume that \ $n \ge 2$.

Then, given \ $[f_i] \in \pi_m (N, *_i), i = 1, 2, \ \ *_1 \ne *_2$, we can identify \ $\Omega^\# (f_1, f_2)$, \  $\Omega_{m -n} (E (f_1, f_2); \widetilde\varphi)$ \ and \ $\Omega_{m -n} (M; \varphi)$ \ with the corresponding groups in the top line of the diagram
$$
   \xymatrix{
     \pi_m(S^n\wedge\Omega(N)^+) \ar[rr]^-{\text{stabilize}}  & &
     \Omega_{m-n}^{\text{fr}}(\Omega N)  \ar[rr]  & &
     \Omega_{m-n}^{\text{fr}}    \\
     & & \pi_m(N) \ar[ull]^{\deg^\#}  \ar[u]_{\widetilde{\deg}}
         \ar[rru]_{\deg}
    .}
   \tag{1.14}
$$
(This is possible since the loop space \ $\Omega N$ \ occurs as a typical fiber of the natural projection \ $p : E (f_1, f_2) \to S^m$, \ cf.\ \cite{K 3}, \S\ 7 and \cite{K 5}).

Furthermore, after deforming the maps \ $f_1$ \ and \ $f_2$ \ until they are constant on opposite halfspheres in \ $S^n$, \ we see that
$$
\widetilde\omega (f_1, f_2) \ = \ \widetilde\omega (f_1, *_2) \ + \widetilde\omega (*_1, f_2),
$$
and similarly for \ $\omega^\#$ \ and \ $\omega$.

Thus it suffices to study the degree maps in diagram 1.14. They turn out to be group homomorphisms which commute with the indicated natural forgetful homomorphisms.

It can be shown (cf. \cite{K 5}) that \ $\deg^\# (f)$ \ is (a strong version of) the Hopf-Ganea invariant of \ $[f] \in \pi_m (N)$ \ (w.r.\ to the attaching map of a top cell in \ $N$, \ compare \cite{CLOT}, 6.7), while \ $\widetilde\deg (f)$ \ is closely related to (weaker) stabilized Hopf-James invariants (\cite{K 3}, 1.14).
\endexample

\subhead Special case: $M = S^m, N = S^n, n \ge 2$\endsubhead  \ Here \ $\deg^\#$ \ is injective and we see that
$$
N (f, *) \le N^\# (f, *) = MCC (f, *) = \left\{ \aligned
0 \ \ \ &\text{if} \ f \ \text{is nulhomotopic} \\
1 \ \ \ &\text{otherwise} \ \ . \endaligned \right.
$$

\flushpar for all maps \ $f : S^m \to S^n$. There are many dimension combinations \ $(m, n)$ \ where the equality \ $N (f, *) = N^\# (f, *)$ \ is also valid for all \ $f$ \ or, equivalently, where \ $\widetilde\deg$ \ is injective (compare e.g.\ our remark 4.10 below or \cite{K 3}, 1.16). However, if \ $n \ne 1, 3, 7$ \ is odd and \ $m = 2n -1$, or if e.g.\ $(m, n) = (8, 4), (9, 4), (9, 3), (10, 4), (16, 8), (17, 8), (10 + n, n)$ \ for \ $3 \le n \le 11$, \ or \ $(24, 6)$, \ then there exists a map \ $f : S^m \to S^n$ \ such that \ $0 = N (f, *) < N^\# (f, *) = 1 $ \ \ (compare \cite{K 3}, 1.17).

\medskip
\subhead Very special case: \ $M = S^3, N = S^2$
\endsubhead \ Here
$$
\widetilde\deg \ \ : \ \ \pi_3 (S^2) \cong \Bbb Z \ \ @>{\ \ \ \ }>> \ \ \Omega^{fr}_1 (\Omega S^2) \ \cong \ \Bbb Z_2 \oplus \Bbb Z
$$
captures the Freudenthal suspension and the classical Hopf invariant of a homotopy class \ $[f]$; therefore \ $\widetilde\deg$ \ is injective (and so is \ $\deg^\#$ \ a fortiori).

On the other hand the invariant \ $\deg (f) \in \Omega^{fr}_1 \cong \Bbb Z_2$ \ (which does not involve the path space \ $E (f, *))$ \ retains only the suspension of \ $f$. \  The corresponding homological invariant \ $\mu (\deg (f)) \in H_1 (S^3; \Bbb Z)$ \ vanishes altogether. \hfill $\square$

\bigskip
Finally let us point out that our approach can also be applied fruitfully to study linking phenomena. Consider e.g.\ a link map
$$
f \ = \ f_1 \amalg f_2 \ : \ M_1 \ \amalg \ M_2 \ \longrightarrow \ N \times \Bbb R
$$
(i.e.\ the closed manifolds \ $M_1$ \ and \ $M_2$ \ have disjoint images). Just as in the case of two disjoint closed curves in \ $\Bbb R^3$ \ the degree of linking can be measured to some extend by the geometry of the overcrossing locus: it consists of that part of the coincidence locus  of the projections to \ $N$, \ where \ $f_1$ \ is bigger than \ $f_2$ \ (w.r.\ to the \ $\Bbb R$--coordinate). Here the normal bordism/pathspace approach yields strong unlinking obstructions which, in addition, turn out to distinguish a great number of different link homotopy classes. Moreover it leads to a natural notion of \ {\it Nielsen numbers for link maps} \ (cf.\ \cite{K 4}).

%% file: gehoni02.tex
\vskip5mm
\specialhead 2. \ \ The path space \ $\bold E \bold( \bold{f_1, f_2}\bold)$
\endspecialhead

A crucial feature of our approach to Nielsen theory is the central role played by the space \ $E (f_1, f_2)$. \ It yields the Nielsen decomposition of coincidence sets in a very natural geometric fashion. In the defining equation 1.7 \ \ \ $N^I$ \ denotes the space of all continuous paths \ $\theta \ : \ I := [0, 1] \to N$ \ with the compact - open topology. The starting point/endpoint fibration \ $N^I \to N \times N$ \ pulls back, via the map
$$
(f_1, f_2) \ \ : \ \ M \ \to \ N \times N \ \ ,
\tag 2.1
$$
to yield the Hurewicz fibration
$$
p \ \ : \ \ E (f_1, f_2) \ \to \ M
\tag 2.2
$$
defined by \ $p (x, \theta) = x$. \ Given a coincidence point \ $x_0 \in M$, \ the fiber \ $p^{-1} (\{x_0\})$ \ is just the loop space \ $\Omega (N, y_0)$ \ of paths in \ $N$ \ starting and ending at \ $y_0 = f_1 (x_0) = f_2 (x_0)$; \ let \ $\theta_0$ denote the constant path at \ $y_0$.

\proclaim{Proposition 2.3} \ The sequence of group homomorphisms
\medskip

\flushpar
\includegraphics[width=125mm]{diag1.eps}

%

\medskip
\flushpar is exact. Moreover, the fiber inclusion \ $\incl : \Omega (N, y_0) \to E (f_1, f_2)$ \ induces a bijection of the sets
$$
R (f_1, f_2) = \pi_1 (N, y_0) / \text{Reidemeister equivalence} \ \ @>{\ \ \ \ \ }>> \pi_0 (E (f_1, f_2))
$$
where two classes \ $[\theta], [\theta'] \in \pi_1 (N, y_0) = \pi_0 (\Omega (N, y_0))$ \ are called Reidemeister equivalent if \ $[\theta'] = f_{1*} (\tau)^{-1} \cdot [\theta] \cdot f_{2*} (\tau)$ \ for some \ $\tau \in \pi_1 (M, x_0)$.
\endproclaim

The proof is fairly evident. In fact, we are dealing here essentially with the long exact homotopy sequence of the fibration \ $p$.

%% file: gehoni03.tex
\vskip5mm
\specialhead 3. \ \ Normal bordism
\endspecialhead

In this section we recall some standard facts about a geometric language which seems well suited to describe relevant coincidence phenomena in arbitrary codimensions.

Let \ $X$ \ be a topological space and let \ $\varphi$ \ be a virtual real vector bundle over \ $X$, i.e.\ an ordered pair \ $(\varphi^+, \varphi^-)$ \ of vector bundles written \ $\varphi = \varphi^+ - \varphi^-$.

{\it A singular \ $\varphi$-manifold in \ $X$ \ of dimension \ $q$} \ is a triple \ $(C, g, \overline g)$ \ where

(i) \ \ \ $C$§ \ is a closed smooth $q$-dimensional manifold;

(ii) \ \ $g : C \to X$ \ is a continuous map;

and

(iii) \ $\overline g : TC \oplus g^* (\varphi^+) \ \longrightarrow \ g^* (\varphi^-)$ \ is a {\it stable} \ vector bundle isomorphism (i.e.\ we can first add trivial vector bundles of suitable dimensions on both sides).

Two such triples \ $(C_i, g_i, \overline g_i), \ i = 0,1$, \ are \ {\it bordant} \ if there exists a compact singular \ $(q + 1)$-dimensional \ $\varphi$-manifold \ $(B, b, \overline b)$ \ in \ $X$ \ with boundary \ $\partial B = C_0 \amalg C_1$ \ such that \ $b$ \ and \ $\overline b$, \ when restricted to \ $\partial B$, \ coincide with the corresponding data \ $g_i$ \ and \ $\overline g_i$ \ at \ $C_i, \ i = 0, 1$ \ (via vector fields pointing into \ $B$ \ along \ $C_0$ \ and out of \ $B$ \ along \ $C_1$). The resulting set of bordism classes, with the sum operation given by disjoint unions, is the \ $q^{th}$ \ {\it normal bordism group \ $\Omega_q (X; \varphi)$ \ of \ $X$ \ with coefficients in \ $\varphi$}.

\example{Example 3.1} \ Let \ $G$ \ denote the trivial group or the (special) orthogonal group \ $(S)O(q'), \ \ q' > q + 1$. For any topological space \ $Y$ \ let \ $\varphi^+$ \ be the classifying bundle over \ $BG$, pulled back to \ $X = Y \times BG$ ,\ while \ $\varphi^-$ \ is trivial. Then \ $\Omega_q (X; \varphi)$ \ is the standard (stably) framed, oriented or unoriented \ $q^{th}$ bordism group of \ $Y$ (compare e.g. \cite{CF}, I.4 and 8).
\endexample

For every virtual vector bundle \ $\varphi$ \ over a topological space \ $X$ \ there are wellknown Hurewicz homomorphisms
$$
\mu \ : \ \Omega_q (X; \varphi) \ \longrightarrow \ H_q (X; \widetilde{\Bbb Z}_\varphi ) , \ \ \ q \in \Bbb Z,
\tag 3.2
$$
into singular homology with local integer coefficients \ $\widetilde{\Bbb Z}_\varphi$ \ (which are twisted like the orientation line bundle \ $\xi_\varphi = \xi_{\varphi^+} \otimes \xi_{\varphi^-} $ \ of \ $\varphi$); they map a normal bordism class \ $[C, g, \overline g]$ \ to the image of the fundamental class \ $[C] \in H_q (C; \widetilde{\Bbb Z}_{TC})$ \ by the induced homomorphism \ $g_*$.

In most cases \ $\mu$ \ leads to a big loss of information. However for \ $q \le 4$ \ this loss can often be measured so that explicit calculations of (and in) \ $\Omega_q (X; \varphi)$ \ are possible (in particular so when \ $\varphi$ \ is highly nontrivial), see theorem 9.3 in \cite{K 1}. We obtain for example

\proclaim{Lemma 3.3} \ Assume \ $X$ \ is pathconnected. Then
\vskip2ex
{\rm (i)} \ \  $\Omega_0 (X; \varphi) \ @>{\mu}>{\ \cong \ }> \ H_0 (X; \widetilde{\Bbb Z}_\varphi) = \left\{ \aligned \Bbb Z \ & \text{if} \ w_1 (\varphi) = 0 ; \\ \Bbb Z_2 & \ \text{else \ .}\endaligned \right.$
\medskip

{\rm (ii)} \ The following sequence is exact:

$\qquad \qquad \Omega_2 (X; \varphi) @>{\mu}>> H_2 (X; \widetilde{\Bbb Z}_\varphi) @>{w_2 (\varphi)}>> \Bbb Z_2 @>{\delta_1}>> \Omega_1 (X; \varphi) @>{\mu}>> H_1 (X; \widetilde{\Bbb Z}_\varphi) \to 0$.
\vskip2mm

Here \ $\delta_1 (1)$ \ is represented by the invariantly parallelized unit circle, together with a constant map, and
$$
\aligned
w_1 (\varphi) \ =& \ w_1 (\varphi^+) + w_1 (\varphi^-) \qquad \qquad \qquad \qquad \qquad  \qquad \ \ and \\
w_2 (\varphi) \ =& \ w_2 (\varphi^+) + w_1 (\varphi^+) w_1 (\varphi^-) + w_2 (\varphi^-) + w_1 (\varphi^-)^2
\endaligned
$$
denote Stiefel-Whitney classes of \ $\varphi$.
\endproclaim

The setting of (normal) bordism groups provides also a first rate illustration of the fact that the geometric and differential topology of manifolds on one hand, and homotopy theory on the other hand, are often but two sides of the same coin. Indeed, if \ $\varphi^-$ \ allows a complementary vector bundle \ $\varphi^{- \bot}$ \ (such that \ $\varphi^- \oplus \varphi^{- \bot}$ \ is trivial), then the wellknown Pontryagin-Thom construction allows us to interpret \ $\Omega_q (X; \varphi), \ q \in \Bbb Z$, \ as a (stable) homotopy group of the Thom space of \ $\varphi^+ \oplus \varphi^{- \bot}$ \ which consists of the total space of \ $\varphi^+ \oplus \varphi^{- \bot}$ \ with one point \lq\lq added at infinity\rq\rq\ (compare e.g. \cite{CF}, I, 11 and 12). Thus the methods of algebraic topology offer another (and often very powerful) approach to computing normal bordism groups (cf. e.g. chapter II of \cite{CF}).

\example{Example 3.4} \ The Thom space of the vector bundle \ $\varphi = \Bbb R^k$  \ over a one-point space is the sphere \ $S^k = \Bbb R^k \cup \{ \infty\}$. Hence the framed bordism group \ $\Omega^{fr}_q := \Omega_q (\{ \text{point}\}; \varphi)$ \ is canonically isomorphic to the stable homotopy group \ $\pi^S_q := \underset{k \to \infty}\to{\lim} \pi_{q +k} (S^k)$ \ of spheres. It is computed and listed e.g. in Toda's tables (in chapter XIV of \cite{T}) whenever \ $q \le 19$. \hfill $\square$
\endexample

For further details and references concerning normal bordism see e.g. \cite{D} or \cite{K 1}. 

%% file: gehoni04.tex
\vskip5mm
\specialhead 4. \ \ The invariants
\endspecialhead

In this section we discuss the invariants \ $\widetilde\omega (f_1, f_2)$ \ and \ $N (f_1, f_2)$ \ based on normal bordism, as well as their sharper (nonstabilized) versions \ $\omega^\# (f_1, f_2)$ \ and \ $N^\# (f_1, f_2)$. We refer to \cite{K 3}  for some of the details and proofs (see also \cite{K 5}).

In the special case when the map \ $(f_1, f_2) : M \to N \times N$ \ is smooth and transverse to the diagonal
$$
\Delta \ = \ \{(y, y) \in N \times N | y \in N \}
\tag 4.1
$$
the coincidence set
$$
C \ = \ C (f_1, f_2) = (f_1, f_2)^{-1} (\Delta) \ = \ \{x \in M | f_1 (x) = f_2 (x) \}
\tag 4.2
$$
is a smooth submanifold of \ $M$. It comes with the maps
$$
  \xymatrix{
     &
     E(f_1, f_2) \ar[d]^{p} \\
     C \ar[ur]^{\widetilde g} \ar[r]^{g}&
     M
}
\tag 4.3
$$
defined by \ $g (x) = x$ \ and \ $\widetilde g (x) = (x,$ constant path at $f_1 (x) = f_2 (x)), \ \ x \in C$.
\bigskip

\midinsert
\noindent
\includegraphics[width=125mm]{figure1.eps}
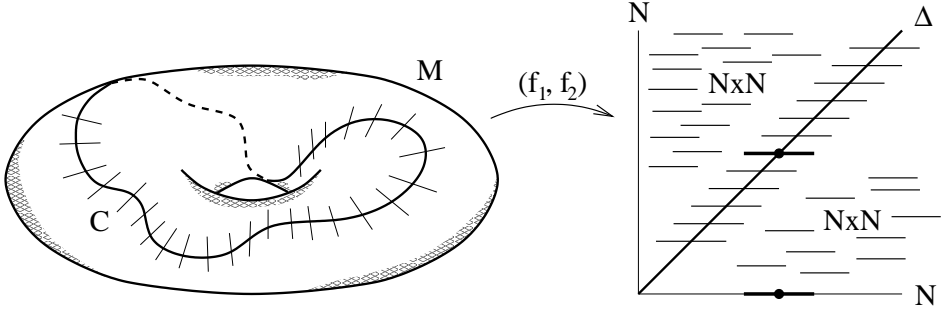
\botcaption{} Figure: A generic coincidence manifold and its normal bundle.
\endcaption
\endinsert

\bigskip
The normal bundle of \ $C$ \ in \ $M$ \ is described by the isomorphism
$$
\nu (C, M) \ \cong  \ (f_1, f_2)^* (\nu (\Delta, N \times N)) \ \cong \ f^*_1 (TN) | C
\tag 4.4
$$
which yields
$$
\overline g : TC \oplus f_1^* (TN) | C  @>{\cong}>> \ TM | C \ .
\tag 4.5
$$

Define
$$
\widetilde\omega (f_1, f_2) \ := \ [C, \widetilde g, \overline g] \in \Omega_{m -n} (E (f_1, f_2); \ \widetilde\varphi)
\tag 4.6
$$
and
$$
\omega (f_1, f_2) \ := \ [C, g, \overline g] = p_* (\widetilde\omega (f_1, f_2) ) \in \Omega_{m -n} (M; \varphi)
\tag 4.7
$$
where
$$
\varphi := f^*_1 (TN) - TM \ \ \text{and} \ \ \widetilde\varphi := p^* (\varphi) .
$$

Invariants with precisely the same properties can be constructed in general. Indeed, apply the preceding procedure to a smooth map \ $(f'_1, f'_2)$ \ which is transverse to \ $\Delta$ \ and approximates \ $(f_1, f_2)$.

Also apply the isomorphism \ $\Omega_* (E (f'_1, f'_2); \widetilde\varphi') \cong \Omega_* (E(f_1, f_2); \widetilde\varphi)$ \ induced by a small homotopy (cf.\ \cite{K 3}, 3.3) to \ $\widetilde\omega (f'_1, f'_2)$ \ in order to obtain \ $\widetilde\omega (f_1, f_2)$\ and similarly \ $\omega (f_1, f_2)$.

Now consider the decomposition
$$
\widetilde\omega (f_1, f_2) = \{ \widetilde\omega_A (f_1, f_2)\} \in \Omega_{m -n} (E(f_1, f_2); \widetilde\varphi) = \bigoplus_A \Omega_{m -n} (A; \widetilde\varphi | A)
\tag 4.8
$$
according to the various pathcomponents \ $A \in \pi_0 (E (f_1, f_2))$ \ of \ $E (f_1, f_2)$.

\definition{Definition 4.9} \ A pathcomponent of \ $E (f_1, f_2)$ \ is called \ {\it essential} \ if the corresponding direct summand of \ $\widetilde\omega (f_1, f_2)$ \ is nontrivial. The \ {\it Nielsen coincidence number} \ $N (f_1, f_2)$ \ is the number of essential pathcomponents \ $A \in \pi_0 (E (f_1, f_2))$.
\enddefinition

Since we assume \ $M$ \ to be compact, \ $N (f_1, f_2)$ \ is a finite integer. It vanishes if and only if \ $\widetilde\omega (f_1, f_2)$ \ does.

\remark{Remark 4.10} \ In the construction above we have neglected an important geometric aspect: \ $C$ \ is much more than just an (abstract) singular manifold with an description of its stable normal bundle. If we keep track (i) of the fact that \ $C$ \ is a smooth \ {\it submanifold} \ in \ $M$, \ and (ii) of the \ {\it nonstabilized} \ isomorphism (4.4), we obtain the sharper invariants \ $\omega^\# (f_1, f_2)$ \  and \ $N^\# (f_1, f_2)$. Note, however, that the bordism set \ $\Omega^\# (f_1, f_2)$ \  in which \ $\omega^\# (f_1, f_2)$ lies has possibly no group structure --- the union of submanifolds may no longer be a submanifold. Also \ $N^\# (f_1, f_2) = 0$ \ if \ $\omega^\# (f_1, f_2) = 0$, \ but the converse may possibly not hold in general --- nulbordisms of individual coincidence components may intersect in \ $M \times I$.

However, in the {\it stable range}\ \ $m \le 2n - 2$ \ \ \ \  $\omega^\# (f_1, f_2)$ \ contains precisely as much information as \ $\widetilde\omega (f_1, f_2)$ \ does, and \ $N^\# (f_1, f_2) = N (f_1, f_2)$.  \hfill $\square$
\endremark
\bigskip

Let us summarize: we have the (successively weaker) invariants \ $\omega^\# (f_1, f_2)$, \ $ \widetilde\omega (f_1, f_2)$, \ \ $\omega (f_1, f_2)$ \ and \ $\mu (\omega (f_1, f_2)) = $ Poincar\'e  dual of the cohomological primary obstruction to deforming \ $f_1$ \ and \ $f_2$ \ away from one another (cf.\ \cite{GJW}, 3.3); they are related by the natural forgetful maps
$$
\Omega^\# (f_1, f_2) @>{\text{stabilize}}>> \Omega_{m -n} (E (f_1, f_2); \widetilde\varphi ) @>{p_*}>> \Omega_{m -n} (M; \varphi) @>{\mu}>> H_{m -n} (M; \widetilde{\Bbb Z}_\varphi)
\tag 4.11
$$
(cf.\ 4.10, 4.3, and 3.2). Only \ $\omega^\# (f_1, f_2)$ \ and \ $\widetilde\omega (f_1, f_2)$ \ involve the path space \  $E (f_1, f_2)$, thus allowing the definition of the Nielsen numbers \ $N^\# (f_1, f_2)$ \ and \ $N (f_1, f_2)$.

\example{Example 4.12: the classical dimension setting  m = n} \ Here the coincidence set
$$
C (f_1, f_2) \ \ = \ \ \coprod_{A \in \pi_0 (E (f_1, f_2))} \ \ \widetilde g^{- 1} (A)
$$

\flushpar consists generically of isolated points (in this very special situation the stabilizing map and the Hurewicz homomorphism \ $\mu$  \ in 4.11 lead to no significant loss of information).

In our approach each Nielsen class is expressed as an inverse image of some pathcomponent \ $A$ \ of \ $E (f_1, f_2)$ \ (compare proposition 2.3). The corresponding index
$$
\widetilde\omega_A (f_1, f_2) \in \Omega_0 (A; \widetilde\varphi|A) \cong \left\{
\aligned &\Bbb Z \ \ \ \ \text{if} \ \omega_1 (\widetilde\varphi | A) = 0 \\
&\Bbb Z_2 \ \ \ \text{else}
\endaligned \right.
$$
(cf.\ lemma 3.3) lies in \ $\Bbb Z_2$ \ precisely if \ $w_1 (\widetilde\varphi | A) \ne 0$ \ or, equivalently, if for some (and hence all) \ $x_0 \in \widetilde g^{-1} (A)$ \ there exists a class \ $\alpha \in \pi_1 (M, x_0)$ \ such that \ $f_{1 *} (\alpha) = f_{2*} (\alpha)$ \ but \ $w_1 (M) (\alpha) \ne f_1^* (w_1 (N)) (\alpha)$ \ (cf.\ \cite{K 3}, 5.2; this agrees with the criterion quoted in \cite{BGZ}, p.\ 53, lines 5--6). If \ $\pi_1 (N)$ \ is commutative, then either the indices of all Nielsen classes are integers, or they all lie in \ $\Bbb Z_2$. However, it is easy to construct examples (e.g.\ involving maps from the Klein bottle to the punctured torus) where both types of pathcomponents  \ $A \in \pi_0 (E (f_1, f_2))$ \ occur.

In any case our approach makes it clear from the outset where the indices of Nielsen classes must take their values.

In the setting of fixed point theory (where \ $f_2$ \ is the identity map on \ $M = N$) \ the transition from the \ $\widetilde\omega$- to the $\omega$-invariant (cf.\ 4.6 and 4.7) which forgets the pathspace \ $E (f_1, f_2)$ \ parallels the transition from Nielsen to Lefschetz theory -- with all the loss of information which this entails.
\hfill $\square$
\endexample 

%% file: gehoni05.tex
\vskip5mm
\specialhead{ 5.\ \  Finiteness conditions for the minimum number \ $\bold M \bold C \bold( \bold f_{\bold 1}, \bold f_{\bold 2} \bold)$}
\endspecialhead

Consider the following possible conditions concerning the invariants defined in 1.2, 4.6, and 4.7:

(C 1) \ \ \ $MC (f_1, f_2) \le 1$ \  ;

(C 2) \ \ \ $MC (f_1, f_2)$ \ is finite \ ;

(C 3) \ \ \ $\widetilde\omega (f_1, f_2)$ \ lies in the image of the homomorphism
$$
i_{E*} \ := \ \bigoplus_A \ i_{A*} \ : \ \bigoplus_A \ \Omega^{fr}_{m -n} \ @>{ \ \ \ \ }>> \ \bigoplus_A \ \Omega_{m -n} (A; \widetilde\varphi|) \ \ = \ \Omega_{m -n} (E (f_1, f_2); \widetilde\varphi )
$$
where direct summation is taken over all \ $A \in \pi_0 (E (f_1, f_2))$ \ and \ $i_{A*}$ \ is induced by the inclusion of a point \ $z_A$ \ into the pathcomponent \ $A$ \ (and by a local orientation of \ $\widetilde\varphi$ \ at \ $z_A$);

(C 4) \ \ \ $\omega (f_1, f_2)$ \ lies in the image of a similarly defined homomorphism
$$
i_* \ : \ \Omega^{fr}_{m -n} \ @>{\ \ \ \ \ }>> \ \Omega_{m -n} (M; \varphi) \ .
$$

\proclaim{Proposition 5.1} \ Each of the first three conditions implies the next one.
\endproclaim

\demo{Proof} \ Assume that the coincidence set \ $C (f_1, f_2)$ \ is finite. If a generic pair \ $(f'_1, f'_2)$ \ approximates \ $(f_1, f_2)$ \ closely enough then each component of \ $C (f'_1, f'_2)$ \ lies in a ball neighbourhood of  some \ $x \in C (f_1, f_2)$; moreover the corresponding paths which occur in the construction of \ $\widetilde\omega (f_1, f_2)$ \ lie entirely in a ball neighbourhood of \ \ $y = f_1 (x) = f_2 (x)$ \ and hence can be contracted into the constant path at $y$. Thus (C 2) implies (C 3). The proposition follows. \hfill $\square$
\enddemo

Our coincidence invariants \ $\widetilde\omega (f_1, f_2)$ \ and \ $\omega (f_1, f_2)$ \ project to the obstructions
$$
\aligned
[\widetilde\omega (f_1, f_2)] \ \ \ &\in \ \ \ \coker (i_{E*}) \\
[\omega (f_1, f_2)] \ \ \ &\in \ \ \ \coker i_*
\endaligned
\tag 5.2
$$
which must vanish if \ $MC (f_1, f_2)$ \ is to be finite.

For \ $m - n = 0$ \ these cokernels are trivial, \ $MC (f_1, f_2)$ \  is actually finite and each integer \ $d \ge 0$ \ can occur as the value of this minimum number for a suitable pair of maps (e.g.\ for selfmaps of degrees \ $d$ \ and \ $0$ \ on \ $S^1$ ).

If \ $m - n = 1$ \ the cokernels in 5.2 are isomorphic --- via the Hurewicz homomorphism \ $\mu$ \ (cf.\ 3.2) --- to \ $H_1 (E (f_1, f_2); \widetilde{\Bbb Z}_{\widetilde\varphi})$ \ and \ $H_1 (M; \widetilde{\Bbb Z}_\varphi)$, \ resp.\ (compare 3.3). In fact \ $\mu$ \ vanishes on the image of \ $i_{E*}$ \ and of \ $i_*$, \ resp., whenever \ $m - n \ge 1$, \ but in general the resulting homomorphisms on the cokernels will not be injective when \ $m - n > 1$ \ (compare \cite{K 3}, 9.3).

\remark{Remark 5.3} \ The finiteness criterion  in 5.1 can be sharpened to yield a nonstabilized version involving \ $\omega^\# (f_1, f_2)$. \hfill $\square$
\endremark
\medskip

For selfcoincidences there is a partial converse of proposition 5.1.

\proclaim{Theorem 5.4} \ If \ $m < 2n - 2$ \ and \ $f_1$ \ is homotopic to \ $f_2$, \ then the four conditions {\rm (C 1) -- (C 4)} \ are equivalent.
\endproclaim

\demo{Proof} \ These conditions are compatible with homotopies of \ $f_1$ \ and \ $f_2$ \ (cf.\ \cite{K 3}, 3.3 and the discussion following 4.4). Hence we may assume that \ $f_1 = f_2 =: f$.

Recall that the selfcoincidence invariant \ $\omega (f, f)$ \ is just the singularity obstruction \ \ $\omega (\underset{\widetilde{\quad}}\to{\Bbb R}, f^* (TN))$ \ to sectioning the vector bundle \ $f^* (TN)$ \ over \ $M$ \ without zeroes (cf.\ \cite{K 2}, theorem 2.2).

 Now assume that \ $m < 2n - 2$ \ and \ $\omega (f, f) = i_* (\omega_0)$ \ for some \ $\omega_0 \in \Omega^{fr}_{m -n} \cong \pi^S_{m -n}$. \ Then there exists a map \ $u^\partial : S^{m -1} \to S^{n -1}$ \ whose (stable) Freudenthal suspension corresponds to \ $\omega_0$. Now consider the trivial bundle \ $f^* (TN) | B^m$ \ over some compact ball \ $B^m$ \ in \ $M$ \ and interpret \ $u^\partial$ \ as a nowhere zero section over \ $\partial B^m = S^{m -1}$.

 We will extend \ $u^\partial$ \ to a section \ $u$ \ of \ $f^* (TN)$ \ over \ {\it all} \ of \ $M$ \ which vanishes only in the centre point of \ $B^m$. \ Over the ball \ $B^m$ \ we use the obvious \lq\lq concentric\rq\rq\ extension. Note, however, that generically the zero set of any extension of \ $u^\partial$ \ over \ $B^m$ \ is a framed manifold which represents \ $\omega_0$. \ Thus a generic extension of \ $u^\partial$ \ to the \ {\it complement} \ $M - \overset{\circ}\to{B}^m$ \ must have a nulbordant manifold of zeroes (representing \ $\omega (f, f) - i_* (\omega_0) = 0)$. \ According to theorem 3.7 in \cite{K 1} these zeroes can be removed altogether.

The resulting section \ $u$ \ of \ $f^* (TN)$ \ allows us to construct a \lq\lq small\rq\rq\ homotopy of \ $f$: \ for every \ $x \in M$ \ just deform \ $f (x)$ \ somewhat in the direction of the tangent vector \ $u (x) \in T_{f (x)} (N)$. \ We obtain a map which has only one coincidence point with \ $f$. \hfill $\square$
\enddemo

%% file: gehoni06.tex
\vskip5mm

\specialhead 6.\ \ The examples of the introduction
\endspecialhead

In view of the selfcoincidence theorem in \cite{K 2} and of our theorem 5.4 the first example is a special case of

\proclaim{Proposition 6.1} \ Let \ $\xi$ \ be an oriented real plane bundle over a closed smooth connected manifold \ $N$ \ and let \ $f \ : \ M := S (\xi) \to N$ \ denote the projection of the corresponding unit circle bundle. Then we have

{\rm (i)} \ \ the selfcoincidence invariant \ $\omega (f, f)$ \ {\rm(cf.\ 4.7)}\ vanishes if and only if the Euler number \ $\chi (N)$ \ of \ $N$ \ lies in \ $e (\xi) (H_2 (N; \Bbb Z)) \subset \Bbb Z$ \ {\rm and} \ $\chi (N)$ \ is even; here \ $e (\xi)$ \ denotes the Euler class of \ $\xi$.

{\rm (ii)} \ the finiteness obstruction \ $\mu (\omega (f, f)) \simeq [\omega (f, f)]$ \ {\rm (cf.\ 5.2} and the subsequent discussion in \S\ {\rm 5)} vanishes if and only if \ $\chi (N) \in e (\xi) (H_2 (N; \Bbb Z))$.
\endproclaim

\demo{Proof} \ We will extend the arguments of section 4 in \cite{K 2}. Consider the commuting diagram
$$
   \xymatrix{
      && {}\save[]-<0cm,3ex>*{~\chi(N)~}="m1"
        \ar@{}[d]|-{\text{\rotatebox{-90}{$\in$}}} \restore
      && {}\save[]-<3em,3ex>*{~\omega(f,f)}="m2"
        \ar@{}[d]|-{\text{\rotatebox{-90}{$\in$}}\qquad} \restore
       \Bbb{Z}_2 \ar[d]^-{\text{incl}_*} \\
      \Omega_2(N;-\xi) \ar[rr]^-{\pitchfork} \ar[d]^-{\mu}&&
      \Omega_0^{\text{fr}}(N)
         \ar[rr]^-{\partial} \ar@{=}[d]&&
      \Omega_1(M;-f^*(\xi)) \ar[d]^-{\mu} \\
      H_2(N;\Bbb{Z}) \ar[rr]^-{e(\xi)} \ar[d]^-{w_2(\xi)} &&
      \Bbb{Z} \ar[rr] &&
      H_1(M;\Bbb{Z}) \ar[d] \\
      \Bbb{Z}_2 && && 0
      \ar@{-->} "m1";"m2"
   }
   \tag{6.2}
$$
Here the vertical exact sequences are as described in lemma 3.3. The horizontal lines are exact Gysin sequences of \ $\xi$ \ in normal bordism and in homology (or, equivalently, oriented bordism), compare \cite{K 1}, 9.20 and 9.4. As was shown in \cite{K 2}, \S\ 4, we have \ $\omega (f, f) = \chi (N) \cdot \partial (1)$. \ Since the Stiefel-Whitney class \ $w_2 (\xi)$ \ is the \ $\mod 2$ \ reduction of \ $e (\xi)$, \ the proposition follows. \ \hfill $\square$
\enddemo

Next let us examine example II. If \ $r \ge 2k \ge 2$ \ then according to the theorem in the introduction of  \cite{K 2} only the \lq\lq highest order part\rq\rq\ 1.11 of the \ {\it complete} \ obstruction \ $\omega (f, f)$ \ (to deforming \ $f$ \ away from itself) survives. Thus the finiteness obstruction \ $[\omega (f, f)]$ \ (cf. 5.2) vanishes. If also \ $k \ge 2$ \ then it follows from theorem 5.4 that \ $MC (f, f)$ \ (and hence also \ $MCC (f, f)$) \ equals \ $0$ \ or \ $1$ \ according to whether \ $\omega (f, f)$ \ vanishes or not. It requires using deep results of homotopy theory and of other branches of algebraic topology to decide which of the two values occur actually, but it can be done at least for \ $k \le 10$ \ (cf.\ \cite{K 2}, \S\ 3). However, the case \ $k = 1$ \ (where we deal with the standard projection from \ $S^{r -1}$ \ to real projective space) is elementary. \hfill $\square$

\medskip
Next we turn to example III where \ $N = (S^1)^n$. \ Use the Lie group structure to replace \ $(f_1, f_2)$ \ by the pair \ $(f, *)$ \ which consists of the quotient \ $f = f_1 \cdot f_2^{- 1}$ \ and of the constant map taking values at the unit element of \ $(S^1)^n$. \ This does not change the coincidence sets and data significantly.

Since each torus \ $(S^1)^k$ \ is a \ $K (\Bbb Z^k, 1)$--space the homotopy class of \ $f$ \ is determined by \ $f_* : H_1 (M; \Bbb Z) \to H_1 ((S^1)^n; \Bbb Z)$. \ Moreover, if the image of \ $f_*$ \ has rank \ $k < n$, then \ $f$ \ factors up to homotopy through the lower dimensional torus \ $(S^1)^k$ \ and hence through \ $(S^1)^n - \{ *\}$. \ On the other hand, if the image of \ $f_*$ \ has rank \ $n$ \ (or, equivalently, the Reidemeister set \ $R (f, *) \cong \coker f_*$ \ is finite) then --- according to 1.10 --- all Reidemeister classes are essential and hence \ $N (f, *) = \# R (f, *)$, \ provided \ $\widetilde\omega (f, *) \ne 0$. \ This holds, in particular, if the invariant
$$
\omega (\coll \circ f, *) \ \in \ \Omega_{m -n} (M; - TM)
$$
which corresponds to the bordism class of the stably coframed manifold \ $C (f, *) = f^{- 1} (\{*\}) \  \subset M$, \ is nontrivial. Here the map
$$
\coll \ : \ N \ = \ (S^1)^n \ @>{ \ \ \ \ \ }>> \ N / (N - \overset{\circ}\to{B}) \ \cong \ S^n
$$
collapses the complement of an open ball to a point. The induced cohomology homomorphisms \ $\coll^*$, \ and \ $(f \circ \coll)^*$, resp., map a generator of \ $H^n (S^n; \Bbb Z)$ \ to the cup product \ $\iota_1 \cdots \iota_n \in H^n ((S^1)^n; \Bbb Z)$, \ and to the Poincar\'e \ dual of \ $\mu (\omega(f, *))$,  resp. (compare 4.11 and \cite{GJW}, 3.3). Our claims concerning example III in the introduction follow now from \S\ 5. \hfill $\square$

\medskip
Finally note that the facts described in example IV follow mainly from the discussion in \cite{K 3} (see 1.14 -- 1.17 as well as sections 7 and 8); the calculation of \ $\Omega^{fr}_1 (\Omega S^2)$ \ can be understood easily with the help of our lemma 3.3. \hfill $\square$

\medskip
Let us put the role of the pathspace \ $E (f_1, f_2)$ \ and its influence on the relative strength of our invariants into perspective (compare diagram 4.11).

In the selfcoincidence situation \ $f_1 = f_2 := f$ \ the fibration \ $p : E (f, f) \to M$ \ allows a \ {\it global} \ section \ $s$ \ by constant paths; therefore \ $\widetilde\omega (f, f) = s_* (\omega (f, f))$ \ is precisely as strong as the (usually much weaker) invariant \ $\omega (f, f)$ \ which does not involve any pathspace data. As examples I and II illustrate, \ $\omega (f, f)$ \ may nevertheless capture decisive and very delicate information (which is also registered to some extend by the Nielsen number \ $N (f, f)$ \ in spite of the fact that it can take only the values \ $0$ \ and \ $1$).

In example III our pathspace approach serves to decompose coincidence sets into Nielsen classes. However, it does not seem to enrich the higher dimensional homotopy theoretical aspects of their data very much (as the torus is aspherical; compare 2.3). Still, all natural numbers can occur here as Nielsen numbers of suitable maps \ $f_1$ \ and \ $f_2$.

In contrast, in example IV the higher dimensional topology of \ $E (f_1, f_2)$ \ turns out to be potentially very rich (e.g.\ when \ $N = S^n, n \ge 2$) \ and able to capture much more than just the decomposition into Nielsen classes. 

%% file: gehonref.tex
\vskip5mm

\subhead Acknowledgement
\endsubhead \newline This work was supported in part by the Deutsche  Forschungsgemeinschaft and AARMS (Canada).

\vskip7mm

\head{References}
\endhead